\documentclass{birkjour}
 \newtheorem{thm}{Theorem}[section]
 \newtheorem{cor}[thm]{Corollary}
 \newtheorem{lem}[thm]{Lemma}
 
 \theoremstyle{definition}
 \newtheorem{defn}[thm]{Definition}
 \theoremstyle{remark}

 \numberwithin{equation}{section}

\begin{document}

%
%
%
%
%
%
%
%
%

\title[Different Volume Computation Methods of Graph Polytopes]{Different Volume Computation Methods of Graph Polytopes}

\author[D. Lee]{Daeseok Lee}
\address{Korea Science Academy of KAIST\\ Pusan 614-822, Korea}
\email{daeseoklee0317@gmail.com}
\author[H.-K. Ju]{Hyeong-Kwan Ju}
\address{%
Department of Mathematics\\
Chonnam National University\\
Gwangju, 500-757 \\
Republic of Korea}
\email{hkju@jnu.ac.kr}


\date{January 1, 2004}

\begin{abstract}
The aim of this work is to introduce several different volume computation methods of the graph polytope associated
with various type of finite simple graphs.  Among them, we obtained the recursive volume formula (RVF)
that is fundamental and most useful to compute the volume of the graph polytope for an arbitrary finite simple graph.
\end{abstract}

\maketitle
\section{Introduction}
B\'{o}na, Ju and Yoshida \cite{[BJY]} enumerated certain weighted graphs with the following conditions :
For a given positive integer $k$, a nonnegative integer $n$ and a simple graph
$G=(VG, EG)$ with $VG=[k]$, where $[k] := \{1, 2, \cdots, k\}$ and $[k]_* :=[k]\cup\{0\}$,
 we consider the set\\
 $$W(n;G) := \{\alpha=(n_1, \cdots, n_k)\in ([n]_*)^k\mid ij\in EG \Rightarrow n_i+n_j
  \leq n\}. $$
 We call $\alpha$ satisfying the conditions in the set given above \textit{(vertex-)
 weighted graph}. In fact, the number of weighted graphs is given by an Ehrhart polynomial of some convex polytope in a unit $k-$hypercube.
 Such a convex polytope is determined uniquely for a given finite simple graph as follows: Let $G = (VG,EG)$  be a simple graph with $VG=[k]$.
 Then the {\bf graph polytope} $P(G)$ associated with the graph $G$ is defined as \\
$$P(G) := \{(x_1, x_2, \cdots, x_k)\in [0,1]^k \mid ij\in EG \Rightarrow x_i + x_j \leq 1 \}$$\\
Our main concerns in this article are the volume computation results of graph polytope associated with many
 types of graphs using several different methods. In order to obtain the volume of the graph polytope we need a certain kernel function $K : [0,1]^2\rightarrow \mathbb{R}$
 defined by the following:\\
$$
 K(s,t) := \left\{
    \begin{array}{cc}
    1, & s+t \leq 1 \\
    0, & \mbox{elsewhere}. \\
    \end{array}
 \right.
$$
Then the volume $vol(G)$ of the polytope $P(G)$ is\\
 $$vol(G)= \int_{Q_n}H(x_1, x_2, \cdots, x_n) dx_1dx_2\cdots dx_n,$$
 where $Q_n = [0, 1]^n$ is the $n$-dimensional unit hypercube,
 $$H(x_1, x_2, \cdots , x_n) = \prod_{ij\in EG}K(x_i, x_j).$$

 From now on, all graphs we mentioned will be finite simple graphs and all polytopes are convex.
 The volume of graph polytope $P(G)$ will be denoted by $vol(G)$ rather than $vol(P(G))$. \\
  Chapter 2 introduces a recursive volume formula for
 the volume of the graph polytope and volume formulae for the graph polytope associated with various types of graphs. Chapter 3 describes the graph joins and
 the corresponding volume formula. Volume of the graph polytope associated with the bipartite graph with certain symmetry is dealt in chapter 4. In chapter 5 we
 use the operator theory to find values for interesting series. In the last chapter we mentioned another way to get the volume of graph polytopes, which  uses  Ehrhart polynomial and series.

\section{Recursive Volume Formula}

One of the key and most fundamental techniques is the {\it recursive volume formula}, or RVF. It is also useful.
Next two lemmas will be used to prove the RVF, and can be shown easily. Applications of RVF are discussed.

\begin{lem}\label{vol1}({\bf polytope partitioning})
Let $P$ be a polytope containing a point $x$  and\\  $d(x,F):=inf\{ d(x,y)| y \in F\}$. Then
\begin{equation}
vol(G)=\sum \frac{d(x,F)vol(F)}{n},
\end{equation}
where the sum runs over all facets of $P(G)$.
\end{lem}

\begin{lem}\label{vol11}
If a graph $G$ has no isolated vertex, then the graph polytope $P(G)$ has no facet of the form $x_i=1.$
In other words, P(G) is only composed of facets of form $x_k=0$ or $x_i+x_j=1$ for $ij\in E$.
\end{lem}

\begin{thm} ({\bf RVF})
Let $G=(VG, EG)$ be a graph with the vertex set $VG=[n]$ and having no isolated vertex. Then
\begin{equation}\label{vol3}
vol(G)=\sum_{i=1}^n \frac{vol(G-i)}{2n},
\end{equation}
where $G-i$ is the graph with the vertex set $[n] \setminus i$ and, accordingly, with the inherited edge set in the original edge set $EG$.\\
\end{thm}

\noindent Proof: Let $x=(\frac{1}{2},\frac{1}{2},\cdots,\frac{1}{2}) \in P(G).$
Then $d(x,F)=\frac{1}{2}$ for the facet $F$ in the hyperplane $x_i=0$ (hence $F=P(G-i)$ in the hyperplane.) and
 $d(x,F)=0$ for all other facets $F$ since $x \in F.$ Note that $vol(G-i)$ is the $(n-1)-$dimensional volume of the graph polytope $P(G-i)$. By the Lemmas \ref{vol1} and \ref{vol11} we have the desired formula (\ref{vol3}).\qed

\begin{cor}[path]\label{PAT}
Let $L_n=([n],EL_n)$ be the path with $EL_n:=\{i(i+1)|i \in [n-1]\}.$ Then
\begin{equation*}
vol(L_n)=\frac{E_{n}}{n!},
\end{equation*}
where $E_n$ is the $n-$th Euler number. Hence its generating function is:
\begin{equation*}
\sum_{n \ge 0} vol(L_n)x^n= secx+tanx.
\end{equation*}
\end{cor}
\noindent Proof: Let $\alpha_n=n!vol(L_{n})$. For convenience, we define $\alpha_0:=0! vol(L_0)=1$. It is obvious that $\alpha_1=1=E_1$.

\begin{equation} \label{PATH}
\begin{split}
vol(L_{n+1}) &= \frac{\sum_{i=0}^n vol(L_i \cup L_{n-i})}{2(n+1)}\\
                &= \frac{\sum_{i=0}^n vol(L_i) vol(L_{n-i})}{2(n+1)}\\
\end{split}
\end{equation} by RVF.

\begin{equation*}
\begin{split}2 \alpha_{n+1} &= 2((n+1)!)vol(L_{n+1})\\
                            &= n! \sum_{i=0}^n vol(L_i)vol(L_{n-i})  ({\text { by the formula(\ref{PATH})}})\\
                            &= \sum_{i=0}^n \binom{n}{i} i!vol(L_i)(n-i)!vol(L_{n-i})\\
                            &= \sum_{i=0}^n \binom{n}{i}\alpha_i \alpha_{n-i} {\text { for }} n \geq 1.\\
\end{split}
\end{equation*}

This says that $\alpha_n$ satisfies the same recurrence relation as $E_n$'s together with $E_0=1$. Hence we get the first conclusion $vol(P(L_n))=\frac{E_n}{n!}.$
Since exponential generating function of the Euler number is $secx+tanx$, the second conclusion follows. \qed

\begin{cor}[cycle]
Let $C_n=([n],E_n)$ be the path with $E_n:=\{i(i+1)|i \in [n](n+1:=1)\}.$ Then
\begin{equation*}
vol(C_n)=\frac{E_{n-1}}{2((n-1)!)}.
\end{equation*}
Hence its generating function is:
\begin{equation*}
\sum_{n \ge 1} vol(C_n)x^n=\frac{x(secx+tanx)}{2}.
\end{equation*}
\end{cor}
\noindent Proof: Removing an edge in the cycle $C_n$ results in a path $L_{n-1}$. Hence RVF and Corollary \ref{PAT} imply both of conclusions. \qed

\begin{cor}[complete graph]
Let $K_n=([n],E_n)$ be the path with $E_n:=\{ij|i,j \in [n]\}.$ Then
\begin{equation*}
vol(K_n)=2^{1-n}.
\end{equation*}
Hence its generating function is:
\begin{equation*}
\sum_{n \ge 1} vol(K_n)x^n=\frac{2x}{2-x}.
\end{equation*}
\end{cor}
\noindent Proof: $vol(K_n)=\frac{n vol(K_{n-1})}{2n} = \frac{vol(K_{n-1})}{2}.$ Since $vol(K_1)=1$, the conclusion follows easily.$ \qed$

\begin{cor}[complete bipartite graph]
Let $K_{s,t}=([s]\cup[t],E_{st})$ be the path with $E_{st}:=\{ij|i \in [s], j \in [t]\}.$ Then
\begin{equation}\label{combip}
vol(K_{s,t})= \frac{1}{\binom{s+t}{s}}.
\end{equation}
\end{cor}
\noindent Proof: Since $vol(K_{s,0})=0=vol(K_{0,t}),$ and
$$vol(K_{s,t})=\frac{s vol(K_{s-1,t})+ t vol(K_{s,t-1})}{2(s+t)},$$
by induction, we have the required formula. \qed

\section{Volume of the graph joins}

Join of graphs $H$ and $K$ is the graph $G=(VG,EG)$ where $VG=VH \cup VK$ and $EG=EH \cup EK \cup \{ xy|x \in VH, y \in VK \}$.
We denote the join of graphs $H$ and $K$ simply by $H+K.$ Obviously, the associative law holds. So we can define
$nG:=G+G+\cdots+G$(add $n$ times). The question is that how we can calculate the volume of $G+H$.
Note that $K_{m+n}=K_m +K_n, \ K_{m,n}=D_m + D_n,$ where $D_n$ is a null graph with $n$ vertices.

\begin{defn}[sliced volume]
Let $G=([n], EG)$ be a graph and $a \in [0,1].$\\
$$vol(G,a):= \int_{[0,a]^n} \prod_{ij \in EG} K(x_i, x_j) dx,$$
\end{defn}
where $K(\cdot,\cdot)$ is defined as:
$$
 K(s,t) := \left\{
    \begin{array}{cc}
    1, & s+t \leq a \\
    0, & \mbox{elsewhere}. \\
    \end{array}
 \right.
$$
\noindent It is obvious that $vol(G))=vol(G,1).$ Next theorem gives us a volume formula
for the graph polytope associated with the joined graph.

\begin{thm}
The volume of the graph polytope associated with the graph $G+H$ is:
\begin{equation}
vol(G+H,c)= \int_0^c \int_0^c \Big(\frac{d}{ds} vol(G,s)\Big) \Big(\frac{d}{dt} vol(H,t)\Big) K(s,t) ds dt.
\end{equation}
\end{thm}

\noindent Proof:
\begin{equation*}
\begin{split}
&vol(G+H,c) \\
&= \int_{[0,c]^{|V(G+H)|}}\prod_{ij \in E(G+H)} K(x_i, x_j) dx \quad (x=(x_1, x_2, \cdots, x_{|V(G+H)|}))\\
&= \int_{[0,c]^{|V(G+H)|}} \Big(\prod_{ij \in E(G)}K(x_i, x_j) \Big)  \Big(\prod_{kl \in E(H)} K(x_k, x_l) \Big)\\
& \qquad K(max\{x_1, x_2, \cdots, x_{|V(G)|} \}, max\{x_{|V(G)|+1}, x_{|V(G)|+2}, \cdots, x_{|V(G+H)|} \})dx   \\
& ( \text{let } s:=max\{x_1, x_2, \cdots, x_{|V(G)|} \}, t:= max\{x_{|V(G)|+1}, x_{|V(G)|+2}, \cdots, x_{|V(G+H)|} \})\\
&= \int_0^c \int_0^c \Big(\frac{d}{ds}  \int_{[0,s]^{|V(G)|}}\prod_{ij \in E(G)} K(x_i, x_j) dx_G  \Big) \\
& \qquad \Big(\frac{d}{dt}  \int_{[0,t]^{|V(H)|}}\prod_{kl \in E(H)} K(x_k, x_l) dx_H \Big) K(s,t)dsdt\\
& (\text{where } dx_G=dx_1 dx_2 \cdots dx_{|V(G)|},\ dx_H=dx_{|V(G)|+1} dx_{|V(G)|+2} \cdots dx_{|V(G+H)|})\\
&= \int_0^c \int_0^c \Big(\frac{d}{ds} vol(G,s)\Big) \Big(\frac{d}{dt} vol(H,t)\Big) K(s,t) ds dt.\\
\end{split}
\end{equation*} \qed
\newpage
\begin{thm}\label{combip2}
Let $\frac{1}{2} \le c \le 1.$ Then the sliced volume of the graph polytope $P(K_{m,n})$ is given by the formula:
$$vol(P(K_{m,n}),c)=c^n (1-c)^m + m \sum_{i=0}^i \binom{n}{i} (-1)^n \frac{c^{m+i}-(1-c)^{m+I}}{m+i}. $$
\end{thm}

\noindent Proof:
\begin{equation*}
\begin{split}
 vol(K_{m,n},c) &= vol(D_m + D_n,c)\\
                   &=  \int_0^c \int_0^c \Big(\frac{d}{ds} vol(D_m,s)\Big) \Big(\frac{d}{dt} vol(D_n,t)\Big) K(s,t) ds dt \\
                   &=  \int_0^c \int_0^c \frac{ds^m}{ds} \frac{dt^n}{dt} K(s,t) ds dt \\
                   &=  \int_0^c  \Big(  \frac{ds^m}{ds}  \int_0^c K(s,t) \frac{dt^n}{dt} dt \Big) ds \\
                   &=  \int_0^c  \Big(  \frac{ds^m}{ds}  \int_0^{min(1-s,c)} \frac{dt^n}{dt} dt \Big) ds \\
                   &=  \int_0^c  \Big(  \frac{ds^m}{ds}  (min(1-s,c))^n \Big) ds \\
                   &=  \int_0^{1-c} c^n \frac{ds^m}{ds} ds +  \int_{1-c}^c (1-s)^n \frac{ds^m}{ds} ds \\
                   &=  c^n (1-c)^m + m \int_{1-c}^c s^{m-1}(1-s)^n ds\\
                   &=  c^n (1-c)^m + m \sum_{i=0}^n \binom{n}{i} (-1)^i \int_{1-c}^c s^{m+i-1} ds \\
                   &=  c^n (1-c)^m + m \sum_{i=0}^n \binom{n}{i} (-1)^i  \frac{c^{m+i}-(1-c)^{m+i}}{m+i}. \\
 \end{split}
\end{equation*} \qed

\begin{cor}
The volume of the graph polytope associated with the complete bipartite graph $K_{m,n}$ is:
\begin{equation}
vol(K_{m,n})=  \sum_{i=0}^n \binom{n}{i} (-1)^i \frac{m}{m+i}.
\end{equation}
\end{cor}
\noindent Proof: Substitute $c=1$ in the Theorem (\ref{combip2}). \qed\\

\noindent The following result is immediate by the third term  from the bottom in the formula $(\ref{combip2})$.
\begin{cor}
$$\sum_{i=0}^n \binom{n}{i}(-1)^i \frac{m}{m+i}=\frac{1}{\binom{m+n}{n}}=\sum_{i=0}^m \binom{m}{i} (-1)^i \frac{n}{n+i}.$$
\end{cor}

According to Rudin(\cite{[PMA]}), the \textit{beta function} $B(r,s)$ is defined as:
$$B(r,s)= \int_0^1 x^{r-1}(1-x)^{s-1} dx=\frac{\Gamma(r) \Gamma(s)}{\Gamma(r+s)},$$
 where $\Gamma(x)$ is the gamma function.
Note that $$vol(K_{m,n})=mB(m,n+1)=nB(m+1,n)$$ from the formula \ref{combip}
and the definition of beta function involving gamma function.
Note that we can give another proof of above equality about beta function from the last integration expression
in the proof.

\begin{thm}[multiple join of a graph]\label{mj}
Let $G$ be a graph, $\frac{1}{2} \le r \le 1$ and $k:=|VG|.$ Then, for any positive integer $n$,
$$ \frac{d}{dr} vol(nG,r)=n(1-r)^{k(n-1)} \frac{d}{dr}vol(G,r).$$
\end{thm}
\noindent Proof:
\begin{equation*}
\begin{split}
 vol(nG,r) &=  \int_0^r \int_0^r \Big(\frac{d}{ds} vol(G,s)\Big) \Big(\frac{d}{dt} vol((n-1)G,t)\Big) K(s,t) ds dt \\
              &=  \int_0^r \frac{d}{ds}(vol(G,s) \Big(\int_0^{min(1-s,r)} \frac{d}{dt} vol((n-1)G,t) dt \Big)  ds \\
              &=  \int_0^r \frac{d}{ds}(vol(G,s) vol((n-1)G,min(1-s,r)) ds \\
              &= (i) + (ii) + (iii)\\
\end{split}
\end{equation*}
\begin{equation*}
\begin{split}
              & \text{where }\\
              &  (i)= \int_0^{1-r} vol((n-1)G,r) \frac{d}{ds} vol(G,s)ds \\
              & \quad \ =vol((n-1)G,r)(1-r)^{k} \ (0 \le 1-r \le \frac{1}{2}) \\
              &  (ii)= \int_{1-r}^{1/2} vol((n-1)G,1-s) \frac{d}{ds} vol(G,s)ds \\
              &  \qquad \ = \int_{1-r}^{1/2} vol((n-1)G,1-s) |VG|s^{|VG|-1} ds \\
              &  (iii)= \int_{1/2}^r vol((n-1)G,1-s) \frac{d}{ds} vol(G,s)ds \\
              &   \qquad \ = \int_{1/2}^r (1-s)^{(n-1)|VG|} \frac{d}{ds} vol(G,s)ds. \\
 \end{split}
\end{equation*}
 Now, we denote $\frac{d}{dr} vol(nG,r)$ simply by $a_n.$ Then
 \begin{equation*}
\begin{split}
 a_n &= \frac{d}{dr} ((i)+(ii)+(iii))\\
     &= (1-r)^k a_{n-1} -k(1-r)^{k-1} vol((n-1)G,r)\\
     & + k(1-r)^{k-1} vol((n-1)G,r) + a_1 (1-r)^{(n-1)k}\\
     &= (1-r)^k a_{n-1} +  (1-r)^{(n-1)k}a_1.\\
  \end{split}
\end{equation*}
 Let $F(x,r)=\sum_{n \ge 1} a_n x^n. $ Then, from the previous recursion formula we get
 $$F(x,r)= \frac{xa_1}{(1-(1-r)^k x)^2} = \sum_{n \ge 1} na_1 (1-r)^{(n-1)k} x^n.$$
 Hence, $$ \frac{d}{dr} vol(nG,r)=n(1-r)^{(n-1)k} \frac{d}{dr}vol(G,r).$$  \qed

\begin{cor}
For the value  $\frac{1}{2} \le r \le 1,$ we have
\begin{equation}\label{compgr}
vol(K_n,r)=2^{1-n}-(1-r)^n \text{ and } vol(K_n)=2^{1-n}.
\end{equation}
\end{cor}
\noindent Proof: Since $K_n=nD_1,$ we have the formula (\ref{compgr}) from the Theorem \ref{mj}. \qed

\begin{cor}
For the value  $\frac{1}{2} \le r \le 1,$ we have
\begin{equation*}
vol(nD_k)=2^{-kn}+n2^{-kn} \frac{1}{\binom{kn}{k}} \sum_{i=0}^{k-1} \binom{kn}{i}.
\end{equation*}
\end{cor}
\noindent Proof:
$$\frac{d}{dt}vol(nD_k,t)=n(1-t)^{k(n-1)} \frac{d}{dt}vol(D_k,t)=knt^{k-1} (1-t)^{k(n-1)}.$$
\begin{equation*}
\begin{split}
vol(nD_k) &= 2^{-kn}+kn \int_{1/2}^1 t^{k-1} (1-t)^{k(n-1)} dt \\
             &= 2^{-kn}+kn \int_{1}^0 \Big( -\frac{1}{2}  \Big) \Big(1 -\frac{t}{2}  \Big)^{k-1} \Big( \frac{t}{2}  \Big)^{k(n-1)} dt \\
             &= 2^{-kn}+kn 2^{-kn} \int_{0}^1 t^{k(n-1)} (2-t)^{k-1} dt\\
             &= 2^{-kn}+kn 2^{-kn} \sum_{i=0}^{k-1} \binom{k-1}{i}   \int_{0}^1 t^{k(n-1)} (1-t)^{i} dt\\
             &= 2^{-kn}+kn 2^{-kn} \sum_{i=0}^{k-1} \binom{k-1}{i}  \frac{(k(n-1))!i!}{(k(n-1)+i+1)!}\\
             &= 2^{-kn}+kn 2^{-kn} \sum_{i=0}^{k-1} \frac{1}{k \binom{kn}{k}}   \binom{kn}{k-1-i}  \\
             &= 2^{-kn}+kn 2^{-kn}  \frac{1}{ \binom{kn}{k}}   \sum_{i=0}^{k-1}    \binom{kn}{i}  \\
\end{split}
\end{equation*} \qed

\section{Volume of bipartite graphs}

\begin{defn}
Let $S_n$ be a set of all permutations of $[n]$ and $\sigma \in S_n.$
$$[0,1]_{\sigma}^n := \{x=(x_1,x_2,\cdots, x_n) \in [0,1]^n | x_{\sigma(1)} \le x_{\sigma(2)} \le \cdots \le x_{\sigma(n)} \}.$$
\end{defn}
\noindent Note that $$[0,1]^n = \bigcup_{\sigma \in S_n} [0,1]_{\sigma}^n $$ and each intersection of two different $[0,1]_{\sigma}^n$ has measure 0 so that
for any measurable function $f$, $$\int_{[0,1]^n} fdx= \sum_{\sigma \in S_n} \int_{[0,1]_{\sigma}^n} fdx$$ and
$$\int_{[0,1]_{\sigma}^n} fdx= \int_0^1 \Big(\int_0^{x_{\sigma(n)}} \Big(\int_0^{x_{\sigma(n-1)}} \cdots \Big(\int_0^{x_{\sigma(2)}} f dx_{\sigma (1)} \Big) \cdots dx_{\sigma (n-2)} \Big) dx_{\sigma (n-1)} \Big)dx_{\sigma (n)}. $$
Let $B=(VB,EB)$ be a bipartite graph with $VB= V_1 \cup V_2, V_1=\{ 1,2,\cdots, n\}, V_2=\{v_1,v_2. \cdots, v_m\}, N_i=\{j \in V_1 | jv_i \in EB \}$.

\begin{thm}\label{bipart}
The volume of the graph polytope associated with the bipartite graph $b$ mentioned above is as follows:
\begin{equation*}
vol(B)= \sum_{\sigma \in S_n} \prod_{i=1}^n \frac{1}{i+ \sum_{j=1}^i \alpha_{j, \sigma}},
\end{equation*}
where
$$\alpha_{i, \sigma}= \text{ number of }(\{ v_k \in V_2 | \sigma(i) \in N_k \} \setminus \cup_{j=1}^{i-1} \{v_k \in V_2 | \sigma(j) \in N_k \}).$$
which means the number of vertices in $V_2$ which the smallest among $\sigma$ values of its neighbors is i.
\end{thm}
\noindent Proof:
\begin{equation*}
\begin{split}
vol(B) & = \int_{[0,1]^n} \prod_{j=1}^m (1-max \{x_i \mid i \in N_j \})dx \\
          & = \sum_{\sigma \in S_n} \int_{[0,1]^n_{\sigma}} \prod_{j=1}^m (1-max \{x_i \mid i \in N_j \})dx \\
          & = \sum_{\sigma \in S_n} \int_{[0,1]^n_{\sigma}} \prod_{j=1}^m (min \{1-x_i \mid i \in N_j \})dx \\
          & = \sum_{\sigma \in S_n} \int_{[0,1]^n_{rev(\sigma)}} \prod_{j=1}^m min \{x_i \mid i \in N_j \}dx \\
          & = \sum_{\sigma \in S_n} \int_{[0,1]^n_{\sigma}} \prod_{j=1}^m min \{x_i \mid i \in N_j \}dx \\
          & = \sum_{\sigma \in S_n} \int_{[0,1]^n_{\sigma}} \prod_{i=1}^n x_{\sigma(i)}^{\alpha_{i,\sigma}} dx \\
          & = \sum_{\sigma \in S_n} \int_0^1(  x_{\sigma(n)}^{\alpha_{n,\sigma}}    \int_0^{x_{\sigma(n)}} [ x_{\sigma(n-1)}^{\alpha_{n-1,\sigma}} \cdots \\
          & \int_0^{x_{\sigma(3)}} (x_{\sigma(2)}^{\alpha_{2,\sigma}}  \int_0^{x_{\sigma(2)}}    ( x_{\sigma(1)}^{\alpha_{1,\sigma}}) dx_{\sigma(1)})   dx_{\sigma(2)} \cdots ] dx_{\sigma(n-1)} )  dx_{\sigma(n)}  \\
          & = \sum_{\sigma \in S_n} \prod_{i=1}^n \frac{1}{i+ \sum_{j=1}^i \alpha_{j, \sigma}}.
\end{split}
\end{equation*} 
where $rev(\sigma)$ represents the opposite order of $\sigma$ in the fifth term.
\qed

An {\it automorphism} of a simple graph $G=(VG,EG)$ is a permutation $\pi$ of $VG$ which has
the property that $uv$ is an edge of $G$ if and only if $\pi(u) \pi(v)$ is an edge of $G$.

\begin{thm}
Assume that the bipartite graph $B$ is symmetric, by the sense that for any permutation $\pi\ $on $V_1$,
there exists a permutation $\sigma$ such that the combination of $\pi$ and $\sigma$
induces an automorphism on $G$. Then,
$$vol(P(B))= n! \prod_{i=1}^n \frac{1}{i+\sum_{j=1}^i \alpha_j}$$
\end{thm}
\noindent Proof:The symmetry of the graph $B$ implies that all $\alpha_{i, \sigma}$s are same
for different $\sigma \in S_n.$ The conclusion follows from the Theorem \ref{bipart}.  \qed

\begin{cor}
Let $B_n$ be the graph that is obtained from the complete bipartite graph $K_{n,n}$ by deleting $n$ disjoint edges.
Then,
 $$vol(B_n)= (1+ \frac{1}{n}) \frac{1}{\binom{2n}{n}}. $$
\end{cor}

\noindent In particular, $vol(B_3)=\frac{1}{15}.$ Note that the bipartite graph $B_3$ is the graph obtained from the 1-skeleton of the cube.

\section{An Application related to the Operator Theory}

We introduce here another interesting fact that uses the linear operator theory to obtain the value of a series described in the theorem below.
One of the results related with the operator theory is the computation
of the $vol(C_n)$, which is referred from  Elkies \cite{[ELK]}.
We will restate the lemma regarding this.

We define $K : [0,1]^2\rightarrow \mathbb{R}$ by the following:\\
$$
 K(s,t) := \left\{
    \begin{array}{cc}
    1, & s+t \leq 1 \\
    0, & \mbox{elsewhere}. \\
    \end{array}
 \right.
$$
Then the volume $vol(G)$ of the polytope $P(G)$ is\\
 $$vol(G)= \int_{Q_n}H(x_1, x_2, \cdots, x_n) dx_1dx_2\cdots dx_n,$$
 where $Q_n = [0, 1]^n$ is the $n$-dimensional unit hypercube,\\

 $$H(x_1, x_2, \cdots , x_n) = \prod_{ij\in E}K(x_i, x_j).$$
 We are interested in the computation of the volume of the polytope $P(G)$ for a given simple graph $G=(V, E)$.\\

We define $\mathcal{K}_n$ inductively as in the following:
$$\mathcal{K}_1(t,s):=K(t,s),$$ and
$$\mathcal{K}_n(t,s):=\int_0^1 \mathcal{K}_1(t,x_1)\mathcal{K}_{n-1}(x_1,s)dx_1 (n \geq 2).$$
Let $T:L^2(0,1) \rightarrow L^2(0,1)$ be a linear operator with kernel $\mathcal{K}_1(\cdot,\cdot)$ on $L^2(0,1)$ defined by:
\begin{equation}
(Tg)(t)=\int_0^1 \mathcal{K}_1(t,s)g(s)ds=\int_0^{1-t} g(s)ds.
\end{equation}
From the definition of $\mathcal{K}_n$ we see that $\mathcal{K}_n(\cdot,\cdot)$ is the kernel function of the linear operator $T^n$ as follows:
\begin{equation}
(T^ng)(t)=\int_0^1 \mathcal{K}_n(t,s)g(s)ds.
\end{equation}

The next lemma gives the spectral decomposition of the linear operator $T$, and also of $T^n$.
Its proof is immediate from the standard linear operator theory.(See also Elkies \cite{[ELK]} or Hutson. et. al. \cite{[HPC]}.)

\begin{lem}
The linear operator $T$ is compact and self-adjoint on $L^2(0,1).$ Its eigen values are $\frac{2}{\pi (4k+1)}(k \in \mathbf{Z})$;
the corresponding eigenfunctions are $cos(\pi (4k+1)/2).$ Moreover, The linear operator $T^n$ is compact and self-adjoint on $L^2(0,1).$ Its eigen values are $(\frac{2}{\pi (4k+1)})^n$ with same corresponding eigenfunctions $cos(\pi (4k+1)/2)$. Each of its eigenvalues for $T$ and $T^n$ is simple.
\end{lem}

Our main goal here is to find the value of certain formula using the operator theory. In fact, it is the $vol(P(C_n))$  which is obtained from the RVF.
By the simple calculations we can get the following formula from the definition of $\mathcal{K}_n$ (see \cite{[ELK]}):
 \begin{equation} \label{TRACE}
vol(C_n)= \int_0^1 \mathcal{K}_n(t,t)dt.
\end{equation}
It turned out that the right hand side of the formula (\ref{TRACE}) is the trace of a trace-class operator $T^n$ over the diagonal, and is equal to
\begin{equation*} \label{TRACE1}
\sum_{k=-\infty}^{\infty} \frac{2^n}{(\pi (4k+1))^n}.
\end{equation*}

Note that this series is absolutely convergent for $n \geq 2$. As a summary we have a theorem:
\begin{thm}
For any integer $n \geq 2$ the following holds:
 \begin{equation*}
 \sum_{k=-\infty}^{\infty} \frac{1} {(4k+1)^n} =\frac{ \pi^n vol(C^n)}{2^n}.
 \end{equation*}
\end{thm} \qed

For the case $n=3,$
 \begin{equation*}
 \begin{split}
 1-\frac{1}{3^3}+\frac{1}{5^3}-\frac{1}{7^3}+\cdots + \frac{(-1)^m}{(2m+1)^3}+\cdots &= \frac{\pi^3 vol(C^3)}{8}\\
                                                                                &= \frac{\pi^3vol(K_3)}{8}\\
                                                                                &= \frac{\pi^3 2^{-2}}{8}\\
                                                                                & =\frac{\pi^3}{32},\\
 \end{split}
 \end{equation*}

meanwhile, for the case $n=4,$
 \begin{equation*}
 \begin{split}
 1+\frac{1}{3^4}+\frac{1}{5^4}+\frac{1}{7^4}+\cdots + \frac{1}{(2m+1)^4}+\cdots &= \frac{\pi^4 vol(C^4)}{16}\\
                                                                                &= \frac{\pi^4vol(K_{2,2})}{16}\\
                                                                                &= \frac{\pi^4}{16\binom{2+2}{2}}\\
                                                                                & =\frac{\pi^4}{96}.\\
 \end{split}
 \end{equation*}


\section{Concluding Remarks and Further Problems}

In fact, we have another volume computation method which comes from the Ehrhart polynomial of $P(G)$.
Let $P$ be an integral convex polytope in ${\mathbb{R}}^d$. Then we call $L_P(t)=|tP \cap {\mathbb{Z}}^d|$ {\it the Ehrhart polynomials}
of $P.$ It is known that, for a given $0/1$-polytope $P$,
$$ vol(P)= \lim_{t \rightarrow \infty} \frac{L_P(t)}{t^d}, \text{ where } d=dim(P),$$
or
$$\frac{f(1)}{d!},\text{ where } \sum_{t=0}^{\infty} L_P(t) x^t =\frac{f(x)}{(1-x)^{d+1}}. $$
(Refer \cite{[BS]} or \cite{[RS]} about this.)
If $G$ is a bipartite graph with $n$ vertices, then its graph polytope $P(G)$ is a $0/1$-polytope of dimension $n$.
Hence, we can get the volume $vol(G)$ from the Ehrhart polynomial $L_{P(G)}(t)$, which we can get by using divided difference technique.
(See B\'{o}na et. al. \cite{[BJY]} for details.)



\end{document}